\documentclass[12pt]{article}
\usepackage{amsmath,amsthm,amssymb,amscd}
\usepackage{latexsym}
\newtheorem{thm}{Theorem}[section]
 
 \newtheorem{lem}[thm]{Lemma}
 \newtheorem{prop}[thm]{Proposition}
 \theoremstyle{definition}
 
 \theoremstyle{remark}

 \numberwithin{equation}{section}

\title
{The Cauchy problem for fully nonlinear parabolic systems on manifolds}

\author{ Hong Huang}

\date{}
\begin{document}
\maketitle

\begin{abstract}
 We show the short-time existence and uniqueness of solutions to the Cauchy problem for fully nonlinear systems of arbitrary even order on closed manifolds which are strongly parabolic at the initial values.   The proof  uses a linearization procedure and a fixed-point argument, and the key ingredient is the well known Schauder estimates for linear, strongly parabolic systems due to Solonnikov.

\vspace*{0.4cm}

key words:  fully nonlinear parabolic systems, manifolds, Cauchy problem, Schauder estimates, linearization. \\


\end{abstract}
\maketitle


\section {Introduction}

In this  note we show the short-time existence and uniqueness of solutions to the Cauchy problem for fully nonlinear systems of arbitrary even order on closed manifolds which are strongly parabolic at the initial values.

 Let $(M^n,g)$ be a closed (i.e. compact and without boundary) Riemannian manifold, and $(E, h)$ be a smooth Euclidean (Riemannian) vector bundle of rank $l$ over $M$ equipped with a  connection $\nabla^E$ compatible with $h$ (sometimes we omit the superscript $E$).  Let   $\Gamma(E)$ be the space of smooth sections of $E$, sometimes one also uses $C^\infty (E)$ to denote the same space. Let $r$ be a positive integer. Let $C^r(E)$ be the space of sections $u$ of $E$ which are continuously differentiable up to order $r$, that is, $\nabla^k u$ is a continuous section of the bundle $E\otimes (T^*M)^{\otimes k}$ for any $0\leq k \leq r$, where $(T^*M)^{\otimes k}$ is the $k$-th tensor power of  $T^*M$.

 Given $T>0$, let
$P_t: \text{Dom}(P_t)\cap C^\infty(E) \rightarrow C^\infty (E)$, $t\in [0, T]$,
be a  smooth family of smooth partial differential operators of order $r$, where $\text{Dom}(P_t)$ is an open subset of $C^r(E)$ on which $P_t$ can act. So it also makes sense to say $P_t: \text{Dom}(P_t) \rightarrow C^0(E)$.
For $u \in \text{Dom}(P_t)$, we may write
\begin{equation}
 P_t(u)(x)=F(x, t, u(x), \nabla u(x), \cdot\cdot\cdot, \nabla^{r}u(x))\in E_x,
\end{equation}
 where $x\in M$ and $t \in [0, T]$.
In local coordinates $(x^1, \cdot\cdot\cdot, x^n)$ of $M$ and local frame $\{e_a\}_{a=1}^l$  of $E$, $u=u^ae_a$, and
\begin{equation}
  P_t(u)(x)=F^a(x^1, \cdot\cdot\cdot, x^n, t, u^1(x),\cdot\cdot\cdot,u^l(x), \frac{\partial u^1}{\partial x^1}(x),\cdot\cdot\cdot, (\frac{\partial}{\partial x^n})^{r} u^l(x))e_a(x),
\end{equation}
where each $F^a$ is a smooth ($C^\infty$) function of its arguments (for smooth $P_t$).

Given $u \in \text{Dom}(P_t)$  and $v\in C^\infty(E)$, the linearization of the operator $P_t$ at $u$ in the direction $v$ is
\begin{equation*}
 P_{t*|u}(v)=\frac{\partial}{\partial s}(P_t(u+sv))|_{s=0}=\lim_{s\rightarrow  0}\frac{P_t(u+sv)-P_t(u)}{s}.
\end{equation*}
Given  a  covector $(x, \xi)\in T^*M$, choose $\phi \in C^\infty(M)$ with $\phi(x)=0$ and $d\phi(x)=\xi$.  Given $v_x \in E_x$, choose $v \in C^\infty (E)$ with $v(x)=v_x$. Define the principal symbol $\sigma(P_{t*|u})$ of $P_{t*|u}$ via
\begin{equation*}
\sigma(P_{t*|u})(x,\xi)v_x=\frac{1}{r!}P_{t*|u}(\phi^rv)(x).
\end{equation*}
(Note that the RHS is independent of the choice of $\phi$ and $v$.) Compare for example pp. 62-63 in \cite {P} and pp. 45-46 in \cite {To}; beware that for the principal symbol there is also a slightly different definition in the literature, see for example p. 113 in \cite{LM}.

From now on let the order $r$ of the operator $P_t$ be even.  The operator $P_t$ is strongly elliptic at $u\in \text{Dom}(P_t)$  if  there exists a constant $\lambda >0$ such that
\begin{equation*}
(-1)^{\frac{r}{2}-1} h(\sigma(P_{t*|u})(x,\xi)v_x, v_x) \geq \lambda |\xi|^{r}|v_x|^2
\end{equation*}
for all $(x, \xi)\in T^*M$ and $v_x \in E_x$. (Note that here our convention is different from that in Definition 4.2 on p. 434 in \cite{Ko}; according to our convention the connection Laplacian $\Delta:=\text {tr}_g\nabla^2$ (instead of $-\Delta$) is strongly elliptic.) The above condition  is also called the Legendre-Hadamard condition. Compare (2.5) in \cite{ADN}, Definition 7 in Chapter VII in \cite{LSU} and the second item in Definition 3.36 in \cite{G}.

The above definitions can also be given with weaker smoothness assumption on $P_t$ and $F$.

Let the vector bundle $E$ be as above,  $u_0$ be a given section of $E$, $P_t$ be a family of differential operators of even order $r$,
consider the Cauchy problem
\begin{equation}
\frac{\partial u}{\partial t}(\cdot,t)=P_ t(u(\cdot,t)),  \hspace*{6mm} u(\cdot,0)=u_0,
\end{equation}
for sections $u(\cdot,t)$ of $E$.

\noindent Various special cases of the following result on (1.3) are  known and considered `standard', but I can not find a detailed proof for the general case in the literature.

\begin{thm} \label{thm 1.1} \ \
Let $M$ be a closed (Riemannian) manifold, $E$ be a (Euclidean) vector bundle over $M$, and $P_t: Dom(P_t)  \cap  C^\infty(E) \rightarrow C^\infty(E)$, $t\in [0, T]$,
be a  smooth family of smooth partial differential operators of even order $r$, where $Dom(P_t)$ is an open subset of $C^r(E)$ on which $P_t$ can act. Suppose $u_0\in Dom(P_0) \cap  C^\infty(E)$ and $P_0$ is strongly  elliptic at $u_0$. Then there is $\delta >0$, such that the Cauchy problem (1.3) has  a unique smooth solution defined on $M \times [0,\delta]$.
\end{thm}

 When $P_0$ is strongly  elliptic at $u_0$, we say (1.3) is strongly parabolic at $u_0$. There are extensive researches on short-time existence  and uniqueness for solutions to quasilinear and fully nonlinear parabolic equations/systems in the literature, but most of them focus on equations/systems in domains in Euclidean spaces, see for example, the books \cite {E1}, \cite {LSU}, \cite {F2}, \cite {Li}, \cite {Lu1} and the references therein; the corresponding works on manifolds are much less, we will only mention a few of them. For the  quasilinear, second order case see, for example,  Theorem 4.51 in Aubin \cite {A} and Section 4.4 in Topping \cite {To} (both without proof), Hamilton \cite {H1}, Section 8 in Chapter 15 of Taylor \cite {T2}, and for the  quasilinear, higher order case see  Mantegazza and Martinazzi \cite {MM}.  There  are results in Baker \cite {Ba} and Lamm \cite {L} (see Main Theorem 1 in \cite {Ba} and Theorem 2.4.5 in \cite {L}) for the  fully nonlinear case similar to our Theorem 1.1, but the Main Theorem 1 in Baker \cite {Ba} imposes an extra symmetry condition on the linearized operator (and moreover, it seems that its proof  is incomplete at some places), and Lamm \cite {L} only considers trivial bundles. Also note that there are  results on the existence and uniqueness of solutions to the Cauchy problem for fully nonlinear second order parabolic systems in Section 7.3 of \cite {T1}.
Theorem 1.1 is extended to  transversely parabolic systems on foliated manifolds in Huang \cite {Hu}.

Many geometric evolution equations/systems are only degenerate parabolic, but by using techniques like the so called DeTurck trick one can often convert them to strongly parabolic equations/systems, then one can apply Theorem 1.1 to get  short time solutions.  For example, Ricci flow, mean curvature flow, cross curvature flow and some fourth-order geometric flows are among them.  In this regard note that Bahuaud-Helliwell \cite {BH} establish short-time existence for some higher-order geometric flows.

Theorem 1.1 follows from  Theorem 1.2 below  (with less  smoothness in the  assumption and conclusion)
 through a standard bootstrap argument.  To state Theorem 1.2 we need to introduce some notation. For $0 < \alpha < 1$ let $C^{r+\alpha}(E)$ be the space of $C^{r+\alpha}$-sections of the bundle $E$. For $0< \tau \leq T$, slightly abusing notation, let $E_\tau \rightarrow M \times [0,\tau]$ be the pullback of the bundle $E\rightarrow M$ via the projection $M \times [0,\tau] \rightarrow M$, and   $C^{r+\alpha, (r+\alpha)/r}(E_\tau)$ be the space of $C^{r+\alpha, (r+\alpha)/r}$-sections of the bundle $E_\tau$.  Here to define the  H\"{o}lder spaces for sections of the vector bundle $E$  we  use parallel transport (defined via the connection $\nabla^E$ and the Levi-Civita connection of $(M,g)$), see  Section 2  for the details (compare for example, p. 483 in \cite {N} and  p. 18 in \cite {Ba} (but note that in the definition of $|u|_{2m,1,\alpha;M_\omega}$ on p. 18 in \cite {Ba} one should also use parallel transport, as $\nabla^{2m,1}u$ is a tensor field (instead of a scalar-valued function) on the manifold $M_\omega$)), which is an extension of that in Euclidean spaces (cf. for example \cite {K}, \cite {LSU}, \cite{L},  \cite {Li}, and \cite {Lu1}). (Note that instead of $C^{r+\alpha, (r+\alpha)/r}$ here, Baker \cite {Ba} uses the notation  $C^{r,1,\alpha}$.)

Now let  a (not necessarily $C^\infty$) differential operator $P_t$, $t\in [0, T]$,  be given again by (1.1). Suppose that $\text {Dom} (P_t)$ is an open set of $C^r(E)$ containing the set
\begin{equation*}
  \{ u \in C^r(E) \hspace*{1mm} |   \hspace*{1mm}
   ||u-u_0||_{C^0}+||\nabla u-\nabla u_0||_{C^0}+\cdot\cdot\cdot+||\nabla^r u-\nabla^r  u_0||_{C^0} \leq R_0\}
\end{equation*}
 for some  $u_0 \in C^{r+\alpha}(E)$ ($0< \alpha < 1$) and $R_0>0$, and  the map $F$ is  twice continuously differentiable,  which means that in each local trivialization, each $F^a$ in (1.2) is a $C^2$-function of its arguments.

\begin{thm} \label{thm 1.2} \ \  Given $0< \alpha < 1$, let  $P_t$, $F$  and $u_0\in C^{r+\alpha}(E)$ be as in the preceding paragraph. Assume that $P_0$ is strongly  elliptic at $u_0$. Then there is $\delta >0$, such that the Cauchy problem (1.3) has  a  unique solution $u \in C^{r+\alpha, (r+\alpha)/r}(E_\delta)$ defined on $M \times [0,\delta]$.
\end{thm}

 Theorem 1.2 is somewhat sharper than the corresponding statements in Main Theorem 1  of Baker \cite {Ba} and in Theorem 2.4.5 of Lamm \cite {L} in that here we do not have any lose on the H\"{o}lder exponent in the conclusion. (Note that Theorem 2.5.7 in \cite{Ge}, which deals with the case of a single equation for a scalar-valued function (instead of a system fora vector-valued function), also has some lose on the H\"{o}lder exponent in the conclusion.) Moreover Theorem 2.4.5 of Lamm \cite {L} has stronger assumption on the regularity of the map $F$.
 (Note that  Lamm \cite {L}  considers the Cauchy-Dirichlet problem on compact manifolds with boundary,  our results can also be extended to this case.)

Our proof of Theorem 1.2 uses a linearization procedure and a fixed-point argument, and the key ingredient is the well known Schauder estimates for linear, strongly parabolic systems  due to Solonnikov \cite {So}  (compare Section 10 of Chapter VII in \cite {LSU} and \cite {S17}); of course, we need to adapt Solonnikov's estimates to the manifold case. Similar  method has been used many times in the literature, see in particular  the proof of  Theorem 8.5.4  in Lunardi's book \cite {Lu1} (see also the proof of Theorem 3.2 in \cite{Lu2}),  which deals with fully nonlinear second order parabolic equations in domains in the Euclidean space; compare  Acquistapace-Terreni \cite {AT}. (In contrast,  the proof of  Theorem 2.4.5 in  Lamm \cite {L} uses implicit function theorem instead of fixed point theorem.)
There are many other works on the  Schauder estimates for linear parabolic systems, see
for example  Friedman \cite {F1} \cite {F2},  Giaquinta-Modica \cite {GM},  Lamm \cite {L} (adapting a method from \cite{Si}), and Schlag \cite {Sc}.
 For recent works on Schauder estimates on linear, strongly parabolic systems with less regularity in the coefficients, see for example, \cite {Bo} and
 \cite {DZ}.

Clearly Theorems 1.1 and 1.2 can be extended to the orbifold case.  Moreover, as in \cite{AT} and \cite{Lu2}, we can also show the continuous dependence of the solutions on
the initial values. The details and applications will be given in a future paper.

In Section 2 we first give two  definitions of parabolic Schauder spaces on closed manifolds, and show their equivalence (see Lemma 2.1), then we briefly reformulate  the linear parabolic Schauder theory on closed manifolds,  in particular, we transfer a special case of Solonnikov's Schauder estimates in domains in the Euclidean space to that on closed manifolds (see Theorem 2.3). Using Theorem 2.3 and the method of continuity, we show the existence and uniqueness of solutions to the Cauchy problem for linear, strongly parabolic systems on closed manifolds; see Theorem 2.4. In Section 3 we prove Theorem 1.2 by using the results in Section 2 and following a strategy in Lunardi \cite{Lu1} \cite{Lu2}. Finally in Section 4 we prove Theorem 1.1 by using Theorem 1.2 and a global Schauder estimate in Solonnikov \cite{So}.

\vspace*{0.4cm}

\noindent {\bf Acknowledgements}. {\hspace*{4mm}}   I'm partially supported by NSFC no.12271040 and NSFC no.11171025.

\section{ The linear parabolic Schauder theory }

As in Section 1, we fix a closed Riemannian manifold $(M,g)$, a Euclidean bundle $(E,h)$ over $M$ with a connection $\nabla$ compatible with $h$, and a positive even number $r$.  For $\tau>0$ let $E_\tau$ be defined as in Section 1. Let  $0< \alpha <1$. For a $C^0$-section $u$ of the bundle $E_\tau$, adapting a definition on p.66 in \cite{So} and Definition 1.1.6 on p.3 in \cite{M}, we define its  spatial H\"{o}lder semi-norm (with exponent $\alpha$) to be
\begin{equation*}
[u]_{\alpha; E_\tau}^{\text{space}}= \mathop{\sup}\limits_{\substack{x,y \in M, x\neq y \\ 0 \leq t \leq \tau}} \frac{|u(x,t)-P_{y,x}u(y,t)|}{d(x,y)^\alpha},
\end{equation*}
where $P_{y,x}u(y,t)$ denotes the parallel transport (defined via the connection $\nabla^E$) of the vector  $u(y,t)\in E_y$  along a minimal geodesic $\gamma$ in $M$ from  $y$ to  $x$, and $d(x,y)$ is the distance between the points $x$ and $y$, define  its temporal H\"{o}lder semi-norm (with exponent $\alpha$) to be
\begin{equation*}
[u]_{\alpha; E_\tau}^{\text{time}}=\mathop{\sup}\limits_{\substack{x \in M \\ 0 \leq s, t \leq \tau, s\neq t}} \frac{|u(x,s)-u(x,t)|}{|s-t|^{\alpha}},
\end{equation*}
its H\"{o}lder semi-norm  to be
\begin{equation*}
[u]_{\alpha; E_\tau}=[u]_{\alpha; E_\tau}^{\text{space}} + [u]_{\alpha/r; E_\tau}^{\text{time}},
\end{equation*}
and its $C^{\alpha, \alpha/r}$-norm to be
\begin{equation*}
||u||_{\alpha, \alpha/r;E_\tau}=||u||_{0; E_\tau} + [u]_{\alpha; E_\tau},
\end{equation*}
where
\begin{equation*}
||u||_{0; E_\tau}=\sup\limits_{(x,t) \in M \times [0, \tau]} |u(x,t)|
\end{equation*}
is the $C^0$-norm of $u$.  Let   $C^{\alpha,\alpha/r}(E_\tau)$ be the space of the sections $u$ of the bundle  $E_\tau$ with $||u||_{\alpha,\alpha/r ; E_\tau} < \infty$, which is a Banach space.

For any nonintegral $\beta >0$, let $\lfloor \beta \rfloor$ be its integer part. For a section $u$ of the bundle $E_\tau$, adapting a definition on p. 66 in \cite{So}, let the spatial H\"{o}lder semi-norm
\begin{equation*}
[u]_{\beta; E_\tau}^{\text{space}}=\sum\limits_{rj+k=\lfloor \beta \rfloor} [\partial_t^j \nabla^k u]_{\beta-\lfloor \beta \rfloor; (E\otimes (T^*M)^{\otimes k})_\tau}^{\text{space}},
\end{equation*}
and the temporal H\"{o}lder semi-norm
\begin{equation*}
 [u]_{\beta/r; E_\tau}^{\text{time}}= \sum\limits_{0<\beta-rj-k<r} [\partial_t^j \nabla^k u]_{(\beta-rj-k)/r; (E\otimes (T^*M)^{\otimes k})_\tau}^{\text{time}},
\end{equation*}
 where $\partial_t=\frac{\partial}{\partial t}$, and $(E\otimes (T^*M)^{\otimes k})_\tau$ is the pullback of the bundle $E\otimes (T^*M)^{\otimes k} \rightarrow M$ via the projection $M \times [0,\tau] \rightarrow M$, (note that in the definition of $[\partial_t^j \nabla^k u]_{\beta-\lfloor \beta \rfloor, (E\otimes (T^*M)^{\otimes k})_\tau}^{\text{space}}$ we use the parallel transport determined via the connection on the bundle $E\otimes (T^*M)^{\otimes k}$ induced by the connection $\nabla^E$ and the Levi-Civita connection of $(M,g)$,) let the  H\"{o}lder semi-norm
 \begin{equation*}
[u]_{\beta; E_\tau}=[u]_{\beta; E_\tau}^{\text{space}} + [u]_{\beta/r; E_\tau}^{\text{time}},
\end{equation*}
  and finally define (cf. p. 91 in \cite{So}) the $C^{\beta, \beta/r}$-norm to be
\begin{equation*}
||u||_{\beta,\beta/r ; E_\tau}=\sum\limits_{rj+k \leq \lfloor \beta \rfloor} ||\partial_t^j \nabla^k u||_{0; (E\otimes (T^*M)^{\otimes k})_\tau} +[u]_{\beta; E_\tau}.
\end{equation*}
 Let   $C^{\beta,\beta/r}(E_\tau)$ be the space of the sections $u$ of the bundle  $E_\tau$ with $||u||_{\beta,\beta/r ; E_\tau} < \infty$, which is also a Banach space.
Compare  also pp. 7-8 in \cite{LSU} and Definition 2.5.2 in \cite{Ge}.

  It turns out  that there is another way to define the parabolic H\"{o}lder norms. For any point $p \in M$, choose a geodesically convex normal chart $U$ around $p$ (cf. Theorem 6.17 in \cite{Le}). Choose an orthonormal frame $\{e_a\}_{a=1}^l$ in the fiber $E_p$  over $p$, and parallel transport (defined via the connection $\nabla^E$) each $e_a$ along radial geodesics in $U$ emanating from $p$, thus get a local orthonormal frame field, still denoted by $\{e_a\}_{a=1}^l$. Since $M$ is compact, we can cover $M$ by a finite number of such local normal charts $(U_i, \varphi_i)$, $1\leq i \leq n_0$, and get  a finite number of local trivializations for $E$ as above. Write $u= \sum_{a=1}^l u^a_{(i)}e_a^{(i)}$  in the $i$-th local trivialization, and consider the $C^{\beta,\beta/r}$-norm of the $\mathbb{R}^l$-valued function $(u^a_{(i)}(\varphi_i^{-1}(\cdot), \cdot))_{1 \leq a \leq l}$ defined on $\varphi_i(U_i) \times [0,\tau]$, (using the relevant definition on p. 66 and p. 91 in \cite{So},) then take the maximum over $1 \leq i \leq n_0$ and get a number, which is denoted by $||u||_{\beta,\beta/r; E_\tau}^{'}$. This is an alternative way to define the $C^{\beta,\beta/r}$-norm of $u$; cf. for example Definitions 1.3.2 and 1.3.4 in \cite{M} (compare also pp. 427-428 in \cite{Ko}).   In effect,  when   one estimates the H\"{o}lder norm of a  section of $E$ one can compare directly the components of the section in  a certain local trivialization  at two
  points in   a local normal chart without using the parallel transport defined by the connection $\nabla^E$ and the Levi-Civita connection of $(M,g)$.  This is justified by the following lemma; compare for example  Proposition 3.21 in \cite {Ba}, \cite {Bu}
  and Chapter 2 (in particular, p. 31) in \cite {Sz}.
  \begin{lem} \label{lem 2.1}  Let $\beta > 0$ be nonintegral. In the above setting, the $C^{\beta,\beta/r}$-norms $||\cdot||_{\beta,\beta/r; E_\tau}^{'}$ of  sections of $E_\tau$ defined by using their  components in local trivializations for $E$ (constructed above) over geodesically convex normal coordinate charts as in the Euclidean case
 are uniformly equivalent to the $C^{\beta,\beta/r}$-norms $||\cdot||_{\beta,\beta/r; E_\tau}$  of  sections of $E_\tau$ more intrinsically defined via parallel transport determined by the connection $\nabla^E$ and the Levi-Civita connection of $(M,g)$  and without using local trivializations.  That is, there is a positive constant $C$ depending on $(M,g)$, $(E, h, \nabla)$, $\alpha$, $r$ and $m$ (but independent of $\tau$), such that for any section $u$ of $E_\tau$,
  \begin{equation*}
 C^{-1}||u||_{\beta,\beta/r; E_\tau}^{'} \leq ||u||_{\beta,\beta/r; E_\tau} \leq  C||u||_{\beta,\beta/r; E_\tau}^{'}.
  \end{equation*}

     \end{lem}
     \noindent {\bf Proof}.
         Since $M$ is compact,  the norms of the curvature tensor of  $(M,g)$  (resp. of $\nabla^E$) and each of its covariant derivatives are bounded. By Theorems A and B in \cite{E} we can control the  Christoffel symbols of the Levi-Civita connection of $(M,g)$ and the connection $\nabla^E$  and their derivatives in normal coordinate charts and local trivializations constructed above. (See also Corollary 4.11 in \cite {H2} and Theorem 1.3 in \cite {He}.) The estimates for the  terms involving various $C^0$-norms are easy,  as locally we can compute the covariant derivatives of any order by using the ordinary derivatives and the Christoffel symbols of the Levi-Civita connection of $(M,g)$ and  the connection $\nabla^E$ and their derivatives.  The estimates for the temporal H\"{o}lder semi-norms are also not difficult.

      Now we treat the spatial H\"{o}lder semi-norms. For simplicity we only write down the argument for the case $\lfloor \beta \rfloor =0$, the general case being similar (then below one should replace $\beta$ by $\beta- \lfloor \beta \rfloor$, $E$  by $E\otimes (T^*M)^{\otimes k}$ and $u$ by $\partial_t^j \nabla^k u$ for suitable nonnegative integers $j$ and $k$ with $rj+k=\lfloor \beta \rfloor$). Since the metric space $(M,d)$ is compact, where $d$ is the distance function induced by the Riemannian metric $g$, the open cover $\{U_i\}_{1\leq i \leq n_0}$ of $M$ constructed above has a Lebesgue number $\delta_0>0$.

      Given $x_0 \in M$ and $V_0 \in E_{x_0}$, the parallel transport of  $V_0$  along a minimal geodesic $\gamma$ (parametrized by the arc length $s$)
      emanating from $x_0$ is defined by solving the  homogeneous linear first order ODE
\begin{equation}
\nabla_{\dot{\gamma}(s)}V(s)=0, \hspace*{4mm}  V(0)=V_0
\end{equation}
 for $V(s) \in E_{\gamma (s)}$, the fiber of $E$ over the point $\gamma (s) \in M$.   Since the connection  $\nabla$ is compatible with the metric $h$, we have
\begin{equation}
|V(s)|=|V_0|
\end{equation}
for the solution $V(s)$ to the ODE (2.1).
In a geodesically convex local normal chart  for $M$ containing $x_0$ and local orthonormal frame field $\{e_a\}$ for $E$ constructed as before over the chosen local normal chart, we can rewrite (2.1) as
  \begin{equation}
 \frac{dV^a(s)}{ds}+\Gamma_{ib}^a(\gamma (s)) \frac{d\gamma^i(s)}{ds}V^b(s)=0,  \hspace*{4mm} V^a(0)=V_0^a,
 \end{equation}
where $V^a$ ($1\leq a \leq l$) (resp. $\gamma^i$ ($1\leq i \leq n)$) are the components of $V$ (resp. $\gamma$) in the local frame field $\{e_a\}$ (resp. local coordinates $(x^1, \cdot\cdot\cdot, x^n)$), and $\nabla_{\frac{\partial}{\partial x^i}}e_a=\Gamma_{ia}^be_b$.

From (2.2),(2.3) and the control on $\Gamma_{ib}^a$ we have
 \begin{equation*}
 |V^a(s)|\leq |V_0|,   \hspace*{4mm} \text{and}  \hspace*{4mm} |V^a(s)-V_0^a|\leq C |V_0|s,
  \end{equation*}
   where the `$| \cdot    |$' on the LHS of the above two inequalities means the absolute value of a real number, and $C$ is a constant depending on $(M,g)$ and $(E, h, \nabla)$.

   Fix any section $u$ of the bundle $E_\tau$.
   Given points $x, y \in M$, if $d(x,y) \geq \delta_0$ (the Lebesgue number determined above), the relevant terms in the spatial H\"{o}lder semi-norms of $u$ are easily controlled;  if $0< d(x,y) < \delta_0$, then $x$ and $y$ are contained in the same local normal chart, thus for any $t\in [0,\tau]$, using the above estimates, we have
   \begin{equation*}
 \frac{|u^a(x,t)-(P_{y,x}u(y,t))^a|}{d(x,y)^\beta}\leq \frac{|u^a(x,t)-u^a(y,t)|}{d(x,y)^\beta}+C |u(y,t)|d(x,y)^{1-\beta},
  \end{equation*}
   and
   \begin{equation*}
  \frac{|u^a(x,t)-u^a(y,t)|}{d(x,y)^\beta}\leq \frac{|u^a(x,t)-(P_{y,x}u(y,t))^a|}{d(x,y)^\beta}+C |u(y,t)|d(x,y)^{1-\beta}.
    \end{equation*}

    Now the desired result follows easily. \hfill{$\Box$}

\vspace*{0.4cm}

Below we will omit the superscript in  $||\cdot||_{\beta,\beta/r; E_\tau}^{'}$.

 Let $\Omega$ and  $\Omega'$ be  bounded domains in $\mathbb{R}^n$  with smooth boundaries and with
 $\overline{\Omega'} \subset \Omega$,  $Q=\Omega \times [0,T]$ and $Q'=\Omega' \times [0,T]$, $l$ be a positive integer, and $r$ be a positive even number. Let   $\beta > 0$ be nonintegral,  and $L_t: C^{r+\beta}(\overline{\Omega}, \mathbb{R}^l)
 \rightarrow C^\beta (\overline{\Omega},\mathbb{R}^l)$ ($t\in [0,T]$)
  be a family of linear, strongly elliptic operator of order $r$ with
   $(L_tu)^a=\sum_{|I|\leq r, b \leq l}A_b^{aI}(x,t)\partial_Iu^b$, where  $I=(i_1, \cdot\cdot\cdot,i_j)$ is a multi-index (here $1\leq i_1, \cdot\cdot\cdot, i_j \leq n$, and $|I|=j$; our usage of the multi-index $I$ is somewhat different from the usual convention, but it is convenient for us),  $A^I=(A_b^{aI})_{1 \leq a, b \leq l} \in C^{\beta,\beta/r}(\overline{Q}, \mathbb{R}^{l^2})$,   and $\partial_Iu^b=\frac{\partial ^{|I|}u^b}{\partial x^{i_1} \cdot\cdot\cdot \partial x^{i_j}}$. Let $\lambda >0$ be the infimum of the constants in the Legendre-Hadamard condition for $L_t$ over $t\in [0,T]$, and
 $\Lambda =\sum_{|I|\leq r} ||A^I||_{C^{\beta,\beta/r}(\overline{Q}, \mathbb{R}^{l^2})}$.

We have the following crucial
 Schauder estimates (interior w.r.t. the space) for the Cauchy problem for linear, strongly parabolic systems on a Euclidean domain due to Solonnikov \cite {So}.
\begin{thm} \label{thm 2.2}   Let the notation and assumptions be as above.
 Fix  $1< p \leq \infty$.
  There   is a constant $C$  depending on  $\Omega$, $\lambda$, $\Lambda$, $r$, $\beta$, $p$ and the distance between $\Omega'$ and $\partial \Omega$
 such that given  $u_0 \in C^{r+\beta}(\overline{\Omega}, \mathbb{R}^l)$ and $f \in C^{\beta,\beta/r}(\overline{Q}, \mathbb{R}^l)$, if
 $u \in C^{r+\beta, (r+\beta)/r}(\overline{Q}, \mathbb{R}^l)$ is a solution to the linear, strongly parabolic system
\begin{equation*}
\frac{\partial u}{\partial t}(x,t)=L_ tu(x,t)+f(x,t)  \hspace*{4mm} \text{in}  \hspace*{2mm} Q, \hspace*{6mm} u(x,0)=u_0(x)   \hspace*{4mm} \text{in}  \hspace*{2mm} \Omega,
\end{equation*}
 then
 \begin{equation}
  ||u||_{C^{r+\beta, (r+\beta)/r}(\overline{Q'}, \mathbb{R}^l)} \leq C(||f||_{C^{\beta,\beta/r}(\overline{Q}, \mathbb{R}^l)}+||u_0||_{C^{r+\beta}(\overline{\Omega},\mathbb{R}^l)}+||u||_{L^p(\overline{Q},\mathbb{R}^l)}).
 \end{equation}
    \end{thm}
\noindent {\bf Proof}.   This is a special case of Theorem 4.11 in \cite {So}. (\cite {So} deals with systems  parabolic in a more general sense; see also \cite{S17}.) In particular, note the last paragraph in the proof of Theorem 4.11 there for our case of
`purely interior cylinder $Q'$'.  (Compare the global Schauder estimates in Theorem 4.9 of \cite {So} (see also Theorem 10.1 in Chapter VII of \cite{LSU}),  Theorem 3.1 in Section 1 of Chapter 3 in \cite {E2} and
 Theorem VI.21 in \cite {EZ}.)
 \hfill{$\Box$}

  \vspace*{0.4cm}
 \noindent {\bf Remark}.  There is an interior (w.r.t. space-time)  estimate in Theorem 2.3.10 of Lamm \cite {L},  and a near-bottom estimate similar to (2.4) above in Theorem 2.3.14 of  \cite {L}, where the RHS of the estimate
 contains $||u||_{L^\infty}$
 instead of $||u||_{L^p}$.   See also Theorem 2.3.23 in \cite {L} for the corresponding global Schauder estimate.

\vspace*{0.4cm}

Now let $(M,g)$ be a closed Riemannian manifold, and $(E,h)$ be a Euclidean vector bundle of rank $l$ over $M$ with a connection compatible with $h$. Let $r$ be a positive even number. Let  $\beta > 0$ be nonintegral,   and $L_t: C^{r+\beta}(E) \rightarrow C^\beta (E)$ ($t\in [0,T]$) be a family of linear,
strongly elliptic operators of order $r$.  For a section $u$ of $E$, in terms of the local trivializations constructed above,  we may write $u=u^ae_a$ and $L_tu=(L_tu)^ae_a$, where $(L_tu)^a=\sum_{|I|\leq r, b \leq l}A_b^{aI}(x,t)\partial_Iu^b$. We assume that for each multi-index $I$, the matrix-valued function $A^I=(A_b^{aI})_{1 \leq a, b \leq l}$ is of the class $ C^{\beta,\beta/r}$. Let $\lambda >0$ be the infimum of the constants in the Legendre-Hadamard condition for $L_t$ over $t\in [0,T]$. Let  $\Lambda$  be the maximum of $\sum_{|I|\leq r} ||A^I||_{C^{\beta,\beta/r}}$ over the finitely many local trivializations, where over each local trivialization, for fixed $I$, $||A^I||_{C^{\beta, \beta/r}}$ is defined as before  by viewing the relevant normal coordinate chart as an open subset in  the Euclidean space.

Let $u_0 \in C^{r+\beta}(E)$ and $f \in C^{\beta,\beta/r}(E_T)$, where $E_T$ is defined in Section 1, consider the linear, strongly parabolic system
\begin{equation}
\frac{\partial u}{\partial t}(\cdot,t)=L_ tu(\cdot,t)+f(\cdot,t),  \hspace*{6mm} u(\cdot,0)=u_0,
\end{equation}
on $M \times [0,T]$ for section $u$ of $E_T$.
We have the following Schauder estimate for linear, strongly parabolic systems on closed manifolds, building on the Euclidean case (Theorem 2.2).
\begin{thm} \label{thm 2.3} With the above notation and assumptions, there exists  a constant $C$ depending on  $M$, $E$, $\lambda$, $\Lambda$,  $\beta$ and $T$, such that
if  $u \in C^{r+\beta, (r+\beta)/r}(E_T)$ is  a solution to the linear, strongly parabolic system (2.5) on $M \times [0,T]$, then
\begin{equation}
 ||u||_{C^{r+\beta, (r+\beta)/r}(E_T)} \leq C (||f||_{C^{\beta,\beta/r}(E_T)}+||u_0||_{C^{r+\beta}(E)}).
\end{equation}
 Moreover, $C$ increases w.r.t. $T$.
\end{thm}
\noindent {\bf Proof}   \hspace*{0.2cm}  First as usual  by considering the system
\begin{equation*}
\frac{\partial v}{\partial t}(\cdot,t)=L_ tv(\cdot,t)+L_tu_0(\cdot)+f(\cdot,t),  \hspace*{6mm} v(\cdot,0)=0,
\end{equation*}
satisfying by the section  $v\in E_T$ given by $v(\cdot,t):=u(\cdot,t)-u_0(\cdot)$, and using the estimates
\begin{equation*}
||A^I\partial_Iu_0||_{C^{\beta,\beta/r}} \leq ||A^I||_{C^{\beta,\beta/r}}\cdot ||u_0||_{C^{|I|+\beta}},
 \end{equation*}
 we can reduce the proof of (2.6) to the case of zero initial value. So below we assume $u_0=0$.  Let $u \in C^{r+\beta, (r+\beta)/r}(E_T)$ be a solution to the linear, strongly parabolic system (2.5) with $u(\cdot, 0)=0$. By the inequality (2.4) in Theorem 2.2 (with $p=2$) and Lemma 2.1 we have
  \begin{equation}
   ||u||_{C^{r+\beta, (r+\beta)/r}(E_T)} \leq C(||f||_{C^{\beta,\beta/r}(E_T)}+||u||_{L^2(E_T)}),
  \end{equation}
 where $C$ depends on $M$, $E$, $\lambda$, $\Lambda$ and  $\beta$. (For the transition between the Euclidean
case and the manifold case one can also consult for example, \cite{M}, \cite {Bu} and the proof of Proposition 3.22 in \cite {Ba}.)  Below we'll use the same $C$ to denote various constants different from line to line.

We try to get rid of the term $||u||_{L^2(E_T)}$ on the RHS of (2.7)  by allowing the constant $C$ to depend also on $T$. (Compare the sentence that follows the inequality (66) on p. 1169 in \cite {Sc}.)
  For $s \in [0,T]$, let
   \begin{equation*}
    U(s)=\int_M|u(x,s)|^2d\mu.
        \end{equation*}
Using
 $ u(x, s)=\int_0^s \frac{\partial u}{\partial t}(x,t) dt$
and (2.7) we have
 \begin{equation*}
  U(s)\leq s^2 \text{vol} (M)||\frac{\partial u}{\partial t}||_{C^0(E_s)}^2\leq C (||f||_{C^{\beta,\beta/r}(E_s)}^2+\int_0^sU(t)dt),
      \end{equation*}
where the constant $C$ depends also on $s$, and increases w.r.t. $s$.  By Bellman-Gronwall inequality we get
\begin{equation*}
  U(s)\leq C ||f||_{C^{\beta,\beta/r}(E_s)}^2.
      \end{equation*}
  It follows that
\begin{equation*}
  ||u||_{L^2(E_T)} \leq C ||f||_{C^{\beta,\beta/r}(E_T)},
      \end{equation*}
which, combined with (2.7), implies the desired estimate.
\hfill{$\Box$}

  \vspace*{0.4cm}

   \noindent {\bf Remark}.    There is a similar estimate in Lemma 2.3.25 of Lamm \cite {L} in the case that $E$ is a trivial bundle. But Lamm's proof there uses the parabolic $L^2$-theory
   (Theorem 2.2.1 in \cite {L}) for which  one needs stronger assumption on the regularity of the coefficients of the operator $L_t$ to be able to
   write $L_t$ in  divergence form.

  \vspace*{0.4cm}

Then, based on Theorem 2.3, we have

\begin{thm} \label{thm 2.4}  Let  the coefficients of the linear, strongly elliptic operator $L_t$ ($t \in [0,T]$) be of the class $C^{\beta,\beta/r}$, $f \in C^{\beta,\beta/r}(E_T)$ and $u_0 \in C^{r+\beta}(E)$.
Then there exists a unique solution $u \in C^{r+\beta, (r+\beta)/r}(E_T)$ to the linear, strongly parabolic system (2.5) on $M \times [0,T]$.
\end{thm}
\noindent {\bf Proof}. The uniqueness follows from the Schauder estimate (Theorem 2.3) for the homogeneous, linear, strongly  parabolic system with zero initial value, which is satisfied by the difference of any two solutions in $C^{r+\beta,  (r+\beta)/r}(E_T)$ to the system (2.5).

The proof of existence is by a standard continuity method, compare the proof  of Proposition 3.25 in \cite {Ba} and Theorem 2.4.3 in \cite {L}.  As before we can assume that $u_0=0$.

Given  $0<\beta' < \beta$,  choose a sequence of smooth sections $f_j$ of $E_T$ such that
$$\begin{array}{l}
||f_j||_{C^{\beta,\beta/r}(E_T)} \leq 2 ||f||_{C^{\beta,\beta/r}(E_T)}   \hspace*{4mm}  \text{and} \\
f_j \rightarrow f \hspace*{4mm}  \text{in} \hspace*{4mm} C^{\beta',\beta'/r}(E_T).
\end{array} $$
First we   solve the simpler systems
\begin{equation*}
\frac{\partial v_j}{\partial t}(\cdot,t)=(-1)^{\frac{r}{2}-1}\Delta^{r/2}(v_j(\cdot,t))+f_j(\cdot,t),  \hspace*{6mm} v_j(\cdot,0)=0
\end{equation*}
via parabolic $L^2$-theory as in Sections 2.2 and 2.3 of  Chapter 2 in  Polden \cite {Po} (or Section 7 in Huisken-Polden \cite {HP}) and Section 3.1.2 of Chapter 3 in \cite{Ba} to get a sequence of smooth sections $v_j$ of $E_T$,   where
$ \Delta:=\text {tr}_g\nabla^2=-\nabla^*\nabla$
 is the connection (rough) Laplacian  defined by using $\nabla(=\nabla^E)$ and the Levi-Civita connection of $(M,g)$ (see for example, pp. 9-11 in \cite {C}, and pp.153-155 in \cite {LM}; here $\nabla^*$ is the formal $L^2$-adjoint of $\nabla$). Note that the results in Sections 2.2 and 2.3 of  Chapter 2 in  \cite {Po} can be easily extended to our situation. In particular, G$\mathring{a}$rding's inequality for our operator $(-1)^{\frac{r}{2}}\Delta^{r/2}$ is easier to prove due to its special structure (cf. Lemma 7.7 in \cite {HP} and Lemma 2.2.1 in \cite {Po} for the very strong elliptic operator and scalar-value function case; compare  also Theorem 3.42 in \cite{G} and Theorem 2.1.6 in \cite{L} for the case of vector value functions in Euclidean domains, and Lemma 3.2 in \cite{Ba} (without a proof)  and Theorem 4.2 on p. 435 in \cite{Ko} for the manifold case, both of the latter four for strong elliptic operators). So there is a constant $C$ depending on $(M,g)$ and $(E,h,\nabla)$ such that for any $\psi \in C^\infty (E)$, there holds
\begin{equation}
\int_M h((-1)^{r/2}\Delta^{r/2}\psi,\psi)d\mu \geq \frac{1}{2}||\psi||^2_{W^{r/2,2}(E)}-C||\psi||^2_{L^2(E)}.
\end{equation}
  Here to prove (2.8) we  use the compatibility of the connection $\nabla$ with the fiber metric $h$, integrate by parts (cf. Proposition 8.1 and its proof in Chapter II of \cite {LM}; in particular, here we  use again the assumption that our manifold $M$ is closed),  exchange derivatives, and see that the leading term of the LHS of (2.8) is given by
  \begin{equation*}
  \int_M g^{i_1j_1}\cdot\cdot\cdot g^{i_{r/2}j_{r/2}}h_{ab} \nabla_{i_1\cdot\cdot\cdot i_{r/2}} \psi ^a \cdot \nabla_{j_1\cdot\cdot\cdot j_{r/2}}\psi ^b d\mu,
  \end{equation*}
  while the lower order terms (whose coefficients involve the curvature of $(E, \nabla)$) can be controlled via interpolation.  Now by Theorem 2.3 we have
\begin{equation}
||v_j||_{C^{r+\beta,(r+\beta)/r}(E_T)} \leq C ||f_j||_{C^{\beta,\beta/r}(E_T)} \leq 2C ||f||_{C^{\beta,\beta/r}(E_T)}.
\end{equation}
By Arzela-Ascoli theorem there is a subsequence of $\{v_j\}$ (still denoted by $\{v_j\}$) and a section $v \in C^{r+\beta', (r+\beta')/r}(E_T)$ such that
 \begin{equation}
 v_j \rightarrow v \hspace*{4mm}  \text{in} \hspace*{4mm} C^{r+\beta', (r+\beta')/r}(E_T).
 \end{equation}
 It follows that $v$ solves the equation
 \begin{equation*}
\frac{\partial v}{\partial t}(\cdot,t)=(-1)^{\frac{r}{2}-1}\Delta^{r/2}(v(\cdot,t))+f(\cdot,t),  \hspace*{6mm} v(\cdot,0)=0.
\end{equation*}
  From (2.9) and (2.10) we have  that $v  \in C^{r+\beta, (r+\beta)/r}(E_T)$  and
 \begin{equation*}
||v||_{C^{r+\beta,(r+\beta)/r}(E_T)}  \leq 2C ||f||_{C^{\beta,\beta/r}(E_T)};
\end{equation*}
compare the proof of Lemma 4.1 below.

  Then applying  the method of continuity (see for example Theorem 5.2 in \cite {GT}) to the path of operators
  \begin{equation*}
  (1-\tau)(\frac{\partial}{\partial t}+(-1)^{r/2}\Delta ^{r/2}) +\tau (\frac{\partial}{\partial t}-L_t), \hspace*{4mm} \tau \in [0,1],
  \end{equation*}
  and using the Schauder estimate  in Theorem 2.3 one solves the system (2.5) (with $u_0=0$).
\hfill{$\Box$}

\vspace*{0.4cm}

Consequently we have the following
\begin{thm} \label{thm 2.5}  Let  the linear, strongly elliptic operator $L_t$ ($t \in [0,T]$), $f$ and $u_0$ be smooth ($C^\infty$).
Then there exists a unique smooth solution  to the linear, strongly parabolic system (2.5) on $M \times [0,T]$.
\end{thm}

\section{ Proof of Theorem 1.2}

With the above preparation, the overall strategy of the proof of Theorem 1.2 follows that of  Theorem 8.5.4  in  \cite {Lu1} and Theorem 3.2 in \cite{Lu2}; cf. also \cite{AT}.

We'll need the following result on interpolatory inclusions, which extends Lemma 5.1.1 in \cite{Lu1}.  (The notation for the function spaces appeared in Lemma 3.1 below is a natural extension of that on pp. 175-176 in \cite{Lu1}.)
\begin{lem} \label{lem 3.1} \ \
Let $r$ be a positive even number, $l$ be a positive integer, $0 < \alpha < 1$, $\tau > 0$, $\Omega$ be a bounded domain in $\mathbb{R}^n$ with smooth boundary,  and $Q=\Omega \times [0,\tau]$. Given any $\mathbb{R}^l$-valued function $u=(u^a)_{1\leq a \leq l} \in C^{r+\alpha, 1}(\overline{Q}, \mathbb{R}^l)$, for any integer $0< k < r$, we have
\begin{equation*}
\nabla^ku:=  (\frac{\partial ^{k}u^a}{\partial x^{i_1} \cdot\cdot\cdot \partial x^{i_k}})_{1\leq i_1,\cdot\cdot\cdot,i_k \leq n}^{1\leq a \leq l} \in C^{0, (r-k+\alpha)/r}(\overline{Q}, \mathbb{R}^{n^kl}),
 \end{equation*}
 and there is a constant $C_k=C_k(\Omega, r, l, \alpha)$ independent of $\tau$ and $u$, such that
 \begin{equation*}
 ||\nabla^ku||_{C^{0, (r-k+\alpha)/r}(\overline{Q}, \mathbb{R}^{n^kl})} \leq C_k ||u||_{C^{r+\alpha, 1}(\overline{Q}, \mathbb{R}^l)}.
  \end{equation*}
\end{lem}

\noindent {\bf Proof}. The proof is similar to that of  Lemma 5.1.1 in \cite{Lu1}. Given $u \in C^{r+\alpha, 1}(\overline{Q}, \mathbb{R}^l)$, let $\tilde{u}$ be the map $t \mapsto \tilde{u}(t)=u(t,\cdot)$, $t \in [0,\tau]$. Then $\tilde{u} \in B([0,\tau]; C^{r+\alpha}(\overline{\Omega}, \mathbb{R}^l)) \cap Lip([0,\tau]; C^\alpha(\overline{\Omega}, \mathbb{R}^l))$ (we follow the notation on p.1 and p.3 in \cite{Lu1}).  For any integer $0< k < r$, the space
$C^k(\overline{\Omega},\mathbb{R}^l)$ belongs to the class $J_{(k-\alpha)/r}$ between $C^\alpha(\overline{\Omega},\mathbb{R}^l)$ and $C^{r+\alpha}(\overline{\Omega},\mathbb{R}^l)$  (see Chapter 1 in \cite{Lu1}, Remark 2 on p. 326 in \cite{Tr}, and Lemma 6.37 in \cite{GT}).
By Proposition 1.1.4 (i) in \cite{Lu1},
 $\tilde{u} \in C^{(r-k+\alpha)/r}([0,\tau]; C^k(\overline{\Omega}, \mathbb{R}^l))$, and
\begin{equation*}
 ||\tilde{u}||_{C^{(r-k+\alpha)/r}([0,\tau]; C^k(\overline{\Omega}, \mathbb{R}^l))} \leq c (||\tilde{u}||_{B([0,\tau]; C^{r+\alpha}(\overline{\Omega}, \mathbb{R}^l))}+ ||\tilde{u}||_{Lip([0,\tau]; C^\alpha(\overline{\Omega}, \mathbb{R}^l))}),
  \end{equation*}
where the constant $c$ is independent of $\tau$ and $u$. Now the result follows.
\hfill{$\Box$}

\vspace*{0.4cm}

  Now we prove Theorem 1.2. Given $u_0 \in C^{r+\alpha}(E)$ as in the assumption of Theorem 1.2.
Let  $\min \{1,T\} \geq \delta> 0$ and $R> 0$ be constants to be chosen later (to satisfy three conditions below), and
\begin{equation*}
Y=\{u\in C^{r+\alpha, (r+\alpha)/r}(E_\delta) \hspace*{1mm} |  \hspace*{1mm} u(\cdot,0)=u_0, ||u-u_0||_{C^{r+\alpha, (r+\alpha)/r}(E_\delta)} \leq R\}.
\end{equation*}

Using Lemma 3.1 and a standard procedure via local trivializations for the bundle $E$ (cf. Lemma 2.1 and its proof), we see that there are various positive constants  depending on $M$, $E$ and $r$ but not on $\delta$ or $R$, such that for any $u \in Y$,
$$\begin{array}{l}
 ||u-u_0||_{C^0}+||\nabla u-\nabla u_0||_{C^0}+\cdot\cdot\cdot+||\nabla^r u-\nabla^r  u_0||_{C^0}\\
 \leq C(\delta+\sum_{k=1}^{r-1}C_k\delta^{(r-k+\alpha)/r}+ \delta^{\alpha/r})||u-u_0||_{C^{r+\alpha, (r+\alpha)/r}(E_\delta)} \\
 \leq \widetilde{C}\delta^{\alpha/r}||u-u_0||_{C^{r+\alpha, (r+\alpha)/r}(E_\delta)} \leq    \widetilde{C}\delta^{\alpha/r} R,
\end{array}
$$
where $C_k$ are constants from Lemma 3.1.

The first condition that we want to impose on  $\delta$ and $R$ is
\begin{equation}
\widetilde{C}\delta^{\alpha/r} R \leq \frac{R_0}{2},
\end{equation}
where  $R_0$ is a positive constant given in the paragraph just before the statement of Theorem 1.2 in Section 1. Then  it makes sense to write $F(\cdot, t, u(\cdot,t), \cdot\cdot\cdot, \nabla^{r}u(\cdot,t))$ for  $u \in Y$ and $t \in [0,\delta]$.

As in the proof of Theorem 8.5.4 in \cite {Lu1}, we go to define a map $G: Y \rightarrow Y$ by setting $G(u)=w$ for $u\in Y$, where $w$ is the unique solution in $C^{r+\alpha, (r+ \alpha)/r}(E_\delta)$ to the linear, strongly parabolic system
$$\begin{array}{l}
\frac{\partial w}{\partial t}=P_{0*|u_0}(w(\cdot,t))+F(\cdot, t, u(\cdot,t), \nabla u(\cdot,t), \cdot\cdot\cdot, \nabla^{r}u(\cdot,t))-P_{0*|u_0}(u(\cdot,t)), \\
 w(\cdot,0)=u_0
\end{array}
$$
on $M \times [0,\delta]$ given by Theorem 2.4.

We'll show that for $u, v \in Y$,
\begin{equation}
||G(u)-G(v)||_{C^{r+\alpha, (r+ \alpha)/r}} \leq C(R)\delta^{\alpha/r}||u-v||_{C^{r+\alpha, (r+ \alpha)/r}}.
\end{equation}

The second condition for $\delta$ and $R$ is
\begin{equation}
C(R)\delta^{\alpha/r}\leq \frac{1}{2}.
\end{equation}

Then, using the triangle inequality, (3.2) and (3.3),  we have
\begin{equation}
||G(u)-u_0||_{C^{r+\alpha, (r+ \alpha)/r}}  \leq \frac{R}{2}+||G(u_0)-u_0||_{C^{r+\alpha, (r+ \alpha)/r}}.
\end{equation}

But  $\hat{v}:=G(u_0)-u_0$ satisfies
$$\begin{array}{l}
\frac{\partial \hat{v}}{\partial t}=P_{0*|u_0}(\hat{v}(\cdot,t))+F(\cdot, t, u_0, \nabla u_0, \cdot\cdot\cdot, \nabla^{r}u_0),\\
\hat{v}(\cdot,0)=0.
\end{array}$$

\noindent By Theorem  2.3, there is a constant $C$ independent of $\delta$ (noting that we require $\delta \leq 1$) such that
\begin{equation}
||\hat{v}||_{C^{r+\alpha, (r+ \alpha)/r}}  \leq C ||F(\cdot, \cdot, u_0, \nabla u_0, \cdot\cdot\cdot, \nabla^{r}u_0)||_{C^{\alpha,  \alpha/r}}=:C'.
\end{equation}
Since by assumption $F(\cdot, \cdot, u_0, \nabla u_0, \cdot\cdot\cdot, \nabla^{r}u_0)$ is of the class $C^2$ on $M \times [0,T]$, the constant $C'$ is also independent of $\delta$.

Inequalities (3.4) and (3.5) imply that
\begin{equation}
||G(u)-u_0||_{C^{r+\alpha, (r+ \alpha)/r}}  \leq \frac{R}{2}+C'.
\end{equation}

The third condition for $R$ (and $\delta$) is that  $R$ is suitably large with
\begin{equation}
\frac{R}{2}\geq C'.
\end{equation}

 Now suppose that $R$ and $0 < \delta \leq\min \{1,T\}$ satisfy all the three conditions (3.1), (3.3) and (3.7) (we first choose $R$ to satisfy (3.7), then choose $\delta$ to satisfy (3.1) and (3.3)), and suppose that (3.2) is true. Then $G:Y\rightarrow Y$ is well-defined and is a $\frac{1}{2}$-contraction. So by the contraction mapping principle $G$ has a unique fixed-point, which clearly solves (1.3).

It remains to show (3.2).

 Note that for $u, v \in Y$, $\tilde{w}:=G(u)-G(v)$ satisfies

$$\begin{array}{l}
\frac{\partial \tilde{w}}{\partial t}=P_{0*|u_0}(\tilde{w}(\cdot,t))+F(\cdot, t, u(\cdot,t), \nabla u(\cdot,t), \cdot\cdot\cdot, \nabla^{r}u(\cdot,t))\\
-F(\cdot, t, v(\cdot,t), \nabla v(\cdot,t), \cdot\cdot\cdot, \nabla^{r}v(\cdot,t))-P_{0*|u_0}(u(\cdot,t)-v(\cdot,t)),\\
\tilde{w}(\cdot,0)=0.
\end{array}$$

\noindent By Theorem 2.3 again, there exists a constant $C$ independent of $\delta$ such that
\begin{equation}
||\tilde{w}||_{C^{r+\alpha, (r+ \alpha)/r}} \leq C||\eta||_{C^{\alpha,  \alpha/r}},
\end{equation}
where
$$\begin{array}{l}
\eta(x,t)=F(x,t, u(x,t), \nabla u(x,t), \cdot\cdot\cdot, \nabla^{r}u(x,t))\\
-F(x, t, v(x,t), \nabla v(x,t), \cdot\cdot\cdot, \nabla^{r}v(x,t))-P_{0*|u_0}(u(x,t)-v(x,t)).
\end{array}$$

We'll show
\begin{equation}
||\eta||_{C^{\alpha,\alpha/r}} \leq C(R) \delta^{\alpha/r}||u-v||_{C^{r+\alpha, (r+ \alpha)/r}}.
\end{equation}
Combined with (3.8), it implies (3.2).

By Lemma 2.1, we can estimate the H\"{o}lder norms of $\eta$ and $u-v$ using their components in each  local trivialization  for $E$ (constructed as before) as in the Euclidean case. So below we assume
that $E=\overline{\Omega} \times \mathbb{R}^l$, where $\Omega$ is  a  bounded domain (say an open ball) in $\mathbb{R}^n$ with smooth boundary, $E_\delta=\overline{\Omega}  \times [0, \delta] \times \mathbb{R}^l$, and $u_0 \in C^{r+\alpha}(\overline{\Omega}, \mathbb{R}^l)$. By choosing $\Omega$ to be suitably small, we may assume that the range of $(u_0, \nabla u_0, \cdot\cdot\cdot, \nabla^r u_0)$  is contained in an open ball
$B(\bar{z}, \frac{R_0}{2})$ in $\mathbb{R}^l \times \mathbb{R}^{nl} \times \cdot\cdot\cdot \times \mathbb{R}^{n^rl}$, where $\bar{z} \in \mathbb{R}^l \times \mathbb{R}^{nl} \times \cdot\cdot\cdot \times \mathbb{R}^{n^rl}$, and $F$ is a $C^2$-map from  $\overline{\Omega} \times  [0, \delta] \times B(\bar{z}, R_0)$ to $\mathbb{R}^l$. We'll denote $B(\bar{z}, R_0)$ by $Z$.
We'll write $F=F(x,t,u, q_k^a, \cdot\cdot\cdot, q_{i_1 \cdot\cdot\cdot i_r}^b)$, where  $(q_k^a)_{1\leq k \leq n}^{1\leq a \leq l} \in \mathbb{R}^{nl}$, $\cdot\cdot\cdot$, $(q_{i_1 \cdot\cdot\cdot i_r}^b)_{1\leq i_1,\cdot\cdot\cdot,i_r \leq n}^{1\leq b \leq l} \in \mathbb{R}^{n^rl}$.

For any  $\check{v} \in C^{r}(E)$, we compute the linearization of $P_0$ at $u_0$ in the direction $\check{v}$,

$$\begin{array}{l}
P_{0*|u_0}(\check{v})\\
=F'_{i_1\cdot\cdot\cdot i_r}(x,0, u_0, \nabla u_0, \cdot\cdot\cdot, \nabla^ru_0)\nabla^r_{i_1\cdot\cdot\cdot i_r}\check{v}+\cdot\cdot\cdot \\
 +F'_k(x, 0,u_0, \nabla u_0, \cdot\cdot\cdot, \nabla^ru_0)\nabla_k\check{v}+
F'_u(x,0, u_0, \nabla u_0, \cdot\cdot\cdot, \nabla^ru_0)\check{v}.
\end{array}$$

\noindent Here, by $F'_{i_1\cdot\cdot\cdot i_r}$ we mean the matrix-valued function $(\frac{\partial F^a}{\partial q_{i_1 \cdot\cdot\cdot i_r}^b})_{1\leq a, b \leq l}$,
$\cdot\cdot\cdot$,  $F'_k$  means the matrix-valued function $(\frac{\partial F^a}{\partial q_k^b})_{1\leq a, b \leq l}$,  and   by $\nabla^r_{i_1\cdot\cdot\cdot i_r}\check{v}$ we mean the vector field $(\frac{\partial^r \check{v}^a}{\partial x^{i_1}\cdot\cdot\cdot \partial x^{i_r}})_{1\leq a \leq l}$, $\cdot\cdot\cdot$,  $\nabla_{k}\check{v}$  means the vector field $(\frac{\partial \check{v}^a}{\partial x^k})_{1\leq a \leq l}$.

For  $|\beta|=1$, let  $D^\beta F$ be  any first order derivative of $F$ w.r.t. $z \in Z$. (We'll interpret  $D^0F$ as $F$ itself.)

 Since $F$ is $C^2$,  there exists a constant $K_1$ independent of $\delta$ such that
$$\begin{array}{l}
\sup\{||D^\beta {F}(\cdot,\cdot, z)||_{C^{\alpha, \alpha/r}} \hspace*{1mm} | \hspace*{1mm} z \in Z, |\beta|=0,1\} \leq K_1.
\end{array}
$$

\noindent Similarly  there also exists a constant $K_2$ independent of $\delta$ such that
$$\begin{array}{l}
|D^\beta F(x,t,z_1)-D^\beta F(x,t,z_2)|\leq K_2|z_1-z_2|,\\
\forall (x,t)\in \overline{\Omega}\times [0,\delta], \hspace*{1mm}  z_1, z_2 \in Z, \hspace*{1mm} |\beta|=0,1.
\end{array}
$$

 Then the remaining estimates are very similar to those  in the proof of Theorem 8.5.4 in \cite {Lu1}; cf. also \cite{AT}. But  for completeness we'll reproduce them below.

For $u, v\in Y$ (defined similarly as in the global case), $\sigma \in [0,1]$, $(x,t)\in \overline{\Omega} \times [0, \delta]$, and $u_0$, let

$$\begin{array}{l}
\xi_\sigma(x,t)=\sigma(u(x,t),\nabla u(x,t), \cdot\cdot\cdot, \nabla^{r}u(x,t))\\
+(1-\sigma)(v(x,t),\nabla v(x,t), \cdot\cdot\cdot, \nabla^{r}v(x,t)),\\
\xi_0(x)=(u_0(x),\nabla u_0(x), \cdot\cdot\cdot, \nabla^{r}u_0(x)).\\
\end{array}$$

\noindent Note that $\xi_\sigma(x,t) \in Z$, and there is a constant $C>0$ such that
$$\begin{array}{l}
|\xi_\sigma(x,t))-\xi_\sigma(y,t))| \leq C(||u_0||_{C^{r+\alpha}(E)}+R) |x-y|^\alpha,\\
|\xi_0(x)-\xi_0(y)| \leq C||u_0||_{C^{r+\alpha}(E)}|x-y|^\alpha.
\end{array}$$

We have
$$\begin{array}{l}
F(x,t, u(x,t), \nabla u(x,t), \cdot\cdot\cdot, \nabla^{r}u(x,t))-F(x,t, v(x,t), \nabla v(x,t), \cdot\cdot\cdot, \nabla^{r}v(x,t))\\
=\int_0^1\frac{d}{d\sigma}F(x,t,\xi_\sigma(x,t))d\sigma\\
=\int_0^1F'_{i_1\cdot\cdot\cdot i_r}(x,t,\xi_\sigma(x,t))(\nabla^r_{i_1\cdot\cdot\cdot i_r}u(x,t)-\nabla^r_{i_1\cdot\cdot\cdot i_r}v(x,t))d\sigma+\cdot\cdot\cdot \\
 +\int_0^1F'_k(x,t,\xi_\sigma(x,t))(\nabla_ku(x,t)-\nabla_kv(x,t))d\sigma\\
 +\int_0^1F'_u(x,t,\xi_\sigma(x,t))(u(x,t)-v(x,t))d\sigma,
\end{array}$$
and
$$\begin{array}{l}
\eta(x,t)=\int_0^1(F_{i_1\cdot\cdot\cdot i_r}'(x,t,\xi_\sigma(x,t))-F'_{i_1\cdot\cdot\cdot i_r}(x,0,\xi_0(x)))
\nabla^r_{i_1\cdot\cdot\cdot i_r}(u-v)(x,t)d\sigma\\
+\cdot\cdot\cdot  +\int_0^1(F'_k(x,t,\xi_\sigma(x,t))-F'_k(x,0,\xi_0(x)))\nabla_k(u-v)(x,t)d\sigma\\
+\int_0^1(F'_u(x,t,\xi_\sigma(x,t))-F'_u(x,0,\xi_0(x)))(u-v)(x,t)d\sigma.
\end{array}$$

In order to estimate $|\eta(x,t)-\eta(x,s)|$ for $0\leq s \leq t \leq \delta$, we add and subtract
$$\begin{array}{l}
\int_0^1F'_{i_1\cdot\cdot\cdot i_r}(x,s,\xi_\sigma(x,s))\nabla^r_{i_1\cdot\cdot\cdot i_r}(u-v)(x,t)d\sigma\\
+\cdot\cdot\cdot  +\int_0^1F'_k(x,s,\xi_\sigma(x,s))\nabla_k(u-v)(x,t)d\sigma\\
+\int_0^1F'_u(x,s,\xi_\sigma(x,s))(u-v)(x,t)d\sigma,
\end{array}$$

\noindent and write
$$\begin{array}{l}
\eta(x,t)-\eta(x,s)\\
=\int_0^1(F'_{i_1\cdot\cdot\cdot i_r}(x,t,\xi_\sigma(x,t))-F'_{i_1\cdot\cdot\cdot i_r}(x,s,\xi_\sigma(x,s)))
\nabla^r_{i_1\cdot\cdot\cdot i_r}(u-v)(x,t)d\sigma\\
+\cdot\cdot\cdot  +\int_0^1(F'_k(x,t,\xi_\sigma(x,t))-F'_k(x,s,\xi_\sigma(x,s)))\nabla_k(u-v)(x,t)d\sigma\\
+\int_0^1(F'_u(x,t,\xi_\sigma(x,t))-F'_u(x,s,\xi_\sigma(x,s)))(u-v)(x,t)d\sigma\\
+\int_0^1(F'_{i_1\cdot\cdot\cdot i_r}(x,s,\xi_\sigma(x,s))-F'_{i_1\cdot\cdot\cdot i_r}(x,0,\xi_0(x)))\\
\cdot(\nabla^r_{i_1\cdot\cdot\cdot i_r}(u-v)(x,t)-\nabla^r_{i_1\cdot\cdot\cdot i_r}(u-v)(x,s))d\sigma\\
+\cdot\cdot\cdot  +\int_0^1(F'_k(x,s,\xi_\sigma(x,s))-F'_k(x,0,\xi_0(x)))\\
\cdot(\nabla_k(u-v)(x,t)-\nabla_k(u-v)(x,s))d\sigma\\
+\int_0^1(F'_u(x,s,\xi_\sigma(x,s))-F'_u(x,0,\xi_0(x)))((u-v)(x,t)-(u-v)(x,s))d\sigma.\\
\end{array}$$

Thus we need to estimate
\begin{equation*}
|D^\beta F(x,t,\xi_\sigma(x,t))-D^\beta F(x,s, \xi_\sigma(x,s))|
\end{equation*}
and
 \begin{equation*}
|D^\beta F(x,s,\xi_\sigma(x,s))-D^\beta F(x,0, \xi_0(x))|
\end{equation*}
for $|\beta|=1$.

We have
 $$\begin{array}{l}
 |D^\beta F(x,t,\xi_\sigma(x,t))-D^\beta F(x,s, \xi_\sigma(x,s))|  \\
 \leq |D^\beta F(x,t,\xi_\sigma(x,t))-D^\beta F (x,s, \xi_\sigma(x,t))|\\
 +|D^\beta F(x,s,\xi_\sigma(x,t))-D^\beta F (x,s, \xi_\sigma(x,s))|\\
 \leq (K_1+K_2([\nabla^{r}u]_{C^{0,\frac{\alpha}{r}}}+[\nabla^{r}v]_{C^{0,\frac{\alpha}{r}}}\\
 +\cdot\cdot\cdot+[\nabla u]_{C^{0,\frac{\alpha}{r}}}+[\nabla v]_{C^{0,\frac{\alpha}{r}}}+ [u]_{C^{0,\frac{\alpha}{r}}}+[v]_{C^{0,\frac{\alpha}{r}}}))(t-s)^\frac{\alpha}{r}\\
 \leq C_1(R)(t-s)^\frac{\alpha}{r},
 \end{array}$$
where the notation $[f]_{C^{0,\frac{\alpha}{r}}}$ follows essentially that on p. 175 in \cite{Lu1}, more precisely, for $f\in C(\overline{\Omega} \times [0,\delta], \mathbb{R}^p)$,
\begin{equation*}
[f]_{C^{0,\alpha/r}}=\mathop{\sup}\limits_{\substack{x \in \overline{\Omega} \\ 0 \leq s, t \leq \delta, s\neq t}} \frac{|f(x,s)-f(x,t)|}{|s-t|^{\alpha/r}},
\end{equation*}
 and similarly
 $$\begin{array}{l}
 |D^\beta F(x,s,\xi_\sigma(x,s))-D^\beta F(x,0, \xi_0(x))|\\
 \leq K_1s^\frac{\alpha}{r}+K_2([\nabla^{r}(u-u_0)]_{C^{0,\frac{\alpha}{r}}}+[\nabla^{r}(v-u_0)]_{C^{0,\frac{\alpha}{r}}}
 +\cdot\cdot\cdot\\
 + [\nabla (u-u_0)]_{C^{0,\frac{\alpha}{r}}}+[\nabla(v-u_0)]_{C^{0,\frac{\alpha}{r}}}
 +[u-u_0]_{C^{0,\frac{\alpha}{r}}}+[v-u_0]_{C^{0,\frac{\alpha}{r}}})\delta^\frac{\alpha}{r}\\
 \leq C_2(R)\delta^\frac{\alpha}{r}.
 \end{array}$$

Now
$$\begin{array}{l}
|\eta(x,t)-\eta(x,s)|\\
\leq C_1(R)(t-s)^\frac{\alpha}{r} (\delta^\frac{\alpha}{r}\sum_{i_1,\cdot\cdot\cdot, i_r=1}^n[\nabla^r_{i_1\cdot\cdot\cdot i_r}(u-v)]_{C^{0,\frac{\alpha}{r}}}
+\cdot\cdot\cdot\\
+\delta^\frac{\alpha}{r}\sum_{k=1}^n[\nabla_k(u-v)]_{C^{0,\frac{\alpha}{r}}}+\delta^\frac{\alpha}{r}[u-v]_{C^{0,\frac{\alpha}{r}}})\\
+C_2(R)\delta^\frac{\alpha}{r} ((t-s)^\frac{\alpha}{r}\sum_{i_1,\cdot\cdot\cdot, i_r=1}^n[\nabla^r_{i_1\cdot\cdot\cdot i_r}(u-v)]_{C^{0,\frac{\alpha}{r}}}+\cdot\cdot\cdot\\
+(t-s)^\frac{\alpha}{r}\sum_{k=1}^n[\nabla_k(u-v)]_{C^{0,\frac{\alpha}{r}}}+(t-s)^\frac{\alpha}{r}[u-v]_{C^{0,\frac{\alpha}{r}}})\\
\leq C_3(R)\delta^\frac{\alpha}{r}(t-s)^\frac{\alpha}{r}||u-v||_{C^{r+\alpha,1+\frac{\alpha}{r}}}.
\end{array}$$

\noindent Letting $s=0$, it implies that
\begin{equation*}
||\eta||_{C^0} \leq C_3(R) \delta^{\frac{2\alpha}{r}}||u-v||_{C^{r+\alpha,1+\frac{\alpha}{r}}}.
\end{equation*}

In order to estimate $|\eta(x,t)-\eta(y,t)|$, we add and subtract
$$\begin{array}{l}
\int_0^1(F'_{i_1\cdot\cdot\cdot i_r}(y,t,\xi_\sigma(y,t))-F'_{i_1\cdot\cdot\cdot i_r}(y,0,\xi_0(y)))\nabla^r_{i_1\cdot\cdot\cdot i_r}(u-v)(x,t)d\sigma\\
+\cdot\cdot\cdot  +\int_0^1(F'_k(y,t,\xi_\sigma(y,t))-F'_k(y,0,\xi_0(y)))\nabla_k(u-v)(x,t)d\sigma\\
+\int_0^1(F'_u(y,t,\xi_\sigma(y,t))-F'_u(y,0,\xi_0(y)))(u-v)(x,t)d\sigma,
\end{array}$$

\noindent and write
$$\begin{array}{l}
\eta(x,t)-\eta(y,t)\\
=\int_0^1(F'_{i_1\cdot\cdot\cdot i_r}(x,t,\xi_\sigma(x,t))-F'_{i_1\cdot\cdot\cdot i_r}(y,t,\xi_\sigma(y,t)))
\nabla^r_{i_1\cdot\cdot\cdot i_r}(u-v)(x,t)d\sigma\\
+\cdot\cdot\cdot  +\int_0^1(F'_k(x,t,\xi_\sigma(x,t))-F'_k(y,t,\xi_\sigma(y,t)))\nabla_k(u-v)(x,t)d\sigma\\
+\int_0^1(F'_u(x,t,\xi_\sigma(x,t))-F'_u(y,t,\xi_\sigma(y,t)))(u-v)(x,t)d\sigma\\
+\int_0^1(F'_{i_1\cdot\cdot\cdot i_r}(y,t,\xi_\sigma(y,t))-F'_{i_1\cdot\cdot\cdot i_r}(y,0,\xi_0(y)))\\
\cdot(\nabla^r_{i_1\cdot\cdot\cdot i_r}(u-v)(x,t)-\nabla^r_{i_1\cdot\cdot\cdot i_r}(u-v)(y,t))d\sigma+\cdot\cdot\cdot \\
 +\int_0^1(F'_k(y,t,\xi_\sigma(y,t))-F'_k(y,0,\xi_0(y)))(\nabla_k(u-v)(x,t)-\nabla_k(u-v)(y,t))d\sigma\\
+\int_0^1(F'_u(y,t,\xi_\sigma(y,t))-F'_u(y,0,\xi_0(y)))((u-v)(x,t)-(u-v)(y,t))d\sigma\\
+\int_0^1(F'_{i_1\cdot\cdot\cdot i_r}(y,0,\xi_0(y))-F'_{i_1\cdot\cdot\cdot i_r}(x,0,\xi_0(x)))\nabla^r_{i_1\cdot\cdot\cdot i_r}(u-v)(x,t)d\sigma\\
+\cdot\cdot\cdot  +\int_0^1(F'_k(y,0,\xi_0(y))-F'_k(x,0,\xi_0(x)))\nabla_k(u-v)(x,t)d\sigma\\
+\int_0^1(F'_u(y,0,\xi_0(y))-F'_u(x,0,\xi_0(x)))(u-v)(x,t)d\sigma.\\
\end{array}$$

Then we estimate
$$\begin{array}{l}
|D^\beta F(x,t,\xi_\sigma(x,t))-D^\beta F(y,t, \xi_\sigma(y,t))|
+|D^\beta F(x,0,\xi_0(x))-D^\beta F(y,0, \xi_0(y))|\\
\leq |D^\beta F(x,t,\xi_\sigma(x,t))-D^\beta F(y,t,\xi_\sigma(x,t))|\\
+|D^\beta F(y,t,\xi_\sigma(x,t))-D^\beta F(y,t,\xi_\sigma(y,t))|\\
+|D^\beta F(x,0,\xi_0(x))-D^\beta F(y,0, \xi_0(x))|+|D^\beta F(y,0,\xi_0(x))-D^\beta F(y,0, \xi_0(y))|\\
\leq (K_1+K_2([\nabla^{r}u]_{C^{\alpha,0}}+[\nabla^{r}v]_{C^{\alpha,0}}+\cdot\cdot\cdot
 +[\nabla u]_{C^{\alpha,0}}+[\nabla v]_{C^{\alpha,0}}
+[u]_{C^{\alpha,0}}+[v]_{C^{\alpha,0}})\\
 +K_1+K_2([\nabla^{r}u_0]_{C^{\alpha,0}}+\cdot\cdot\cdot+[\nabla u_0]_{C^{\alpha,0}}+[u_0]_{C^{\alpha,0}}))|x-y|^\alpha\\
 \leq C_4(R)|x-y|^\alpha,
\end{array}$$
\noindent and
$$\begin{array}{l}
|D^\beta F(y,t,\xi_\sigma(y,t))-D^\beta F(y,0, \xi_0(y))|
 \leq C_2(R)\delta^\frac{\alpha}{r}
\end{array}$$
as before.

Now
$$\begin{array}{l}
|\eta(x,t)-\eta(y,t)|\\
\leq C_4(R)(\delta^{\alpha/r}\sum_{i_1,\cdot\cdot\cdot, i_r=1}^n[\nabla^r_{i_1\cdot\cdot\cdot i_r}(u-v)]_{C^{0,\frac{\alpha}{r}}}
+\cdot\cdot\cdot\\
+\delta^{\alpha/r}\sum_{k=1}^n[\nabla_k(u-v)]_{C^{0,{\alpha/r}}}+\delta^{\alpha/r}[u-v]_{C^{0,{\alpha/r}}})|x-y|^\alpha\\
+C_2(R)\delta^{\alpha/r} (\sum_{i_1,\cdot\cdot\cdot, i_r=1}^n[\nabla^r_{i_1\cdot\cdot\cdot i_r}(u-v)]_{C^{\alpha,0}}+\cdot\cdot\cdot\\
+\sum_{k=1}^n[\nabla_k(u-v)]_{C^{\alpha,0}}+[u-v]_{C^{\alpha,0}})|x-y|^\alpha\\
\leq C_5(R)\delta^{\alpha/r}||u-v||_{C^{r+\alpha,1+\frac{\alpha}{r}}}|x-y|^\alpha.
\end{array}$$

Finally combining the above estimates for $\eta$ we get
\begin{equation*}
||\eta||_{C^{\alpha,\frac{\alpha}{r}}} \leq C_6(R) \delta^{\frac{\alpha}{r}}||u-v||_{C^{r+\alpha,1+\frac{\alpha}{r}}},
\end{equation*}
and as said above the short time existence of the solution to (1.3) in $C^{r+\alpha,1+\frac{\alpha}{r}}(E_\delta)$  follows.

The uniqueness of the solution to (1.3) in $C^{r+\alpha, (r+\alpha)/r}(E_\delta)$ follows from the above proof and a `standard' argument as on p. 410 of Lunardi \cite {Lu2}.  We argue by contradiction. Suppose there are two different solutions $u$ and $v$ of the class $C^{r+\alpha, (r+\alpha)/r}(E_\delta)$ to the Cauchy problem (1.3) with the same initial value $u_0 \in C^{r+\alpha}(E)$.  Let
\begin{equation*}
t_0=\sup \{ t \in [0,\delta] \hspace*{1mm} | \hspace*{1mm}  u=v \hspace*{1mm} \text{in}   \hspace*{1mm} M\times [0,t]\}.
\end{equation*}
 Then  $t_0 < \delta$ and $u(x,t_0)=v(x, t_0)$ for any $x \in M$. Set $w_0(x)=u(x,t_0)$ for $x\in M$. Consider the Cauchy problem
 \begin{equation}
\frac{\partial w}{\partial t}(x,t)=P_ t(w(x,t)),  \hspace*{2mm} x\in M,   \hspace*{1mm} t\geq t_0,  \hspace*{6mm} w(x,t_0)=w_0(x), \hspace*{2mm} x\in M.
\end{equation}
    If $R'$ is sufficiently large and $\delta'>0$ is sufficiently small, by the above argument, (3.10) has a unique solution in the space
 \begin{equation*}
Y'=\{w\in C^{r+\alpha, (r+\alpha)/r}(E_{[t_0, t_0+\delta']}) \hspace*{1mm} |  \hspace*{1mm} w(\cdot,t_0)=w_0, ||w-w_0||_{C^{r+\alpha, (r+\alpha)/r}(E_{[t_0, t_0+\delta']})} \leq R'\},
\end{equation*}
  where $E_{[t_0, t_0+\delta']} \rightarrow M \times [t_0, t_0+\delta']$ is the pullback of the bundle $E\rightarrow M$ via the projection $M \times [t_0, t_0+\delta'] \rightarrow M$. We can choose
  \begin{equation*}
R'> \max \{||u-w_0||_{C^{r+\alpha, (r+\alpha)/r}(E_{[t_0, \delta]})}, \hspace*{1mm} ||v-w_0||_{C^{r+\alpha, (r+\alpha)/r}(E_{[t_0, \delta]})}  \}.
\end{equation*}
  Then we get that $u=v$ in $M \times  [t_0, t_0+\delta']$. A contradiction.
  \hfill{$\Box$}

\vspace*{0.4cm}

\noindent {\bf Remark}.  To guarantee the existence of the constants $K_1$ and $K_2$ in the above proof, one only needs that $F$ and $D^\beta F$ ($|\beta|=1$) are locally Lipschitz continuous w.r.t. $z \in Z$ and locally $C^{\alpha, \alpha/r}$ w.r.t. $(x,t)$, uniformly w.r.t. the other variables (cf. Section 8.5.3 in \cite {Lu1}). If one uses the full power of the assumption that $F$ is $C^2$, then the above proof can be slightly simplified, cf. the proof of Theorem 3.2 in \cite {Lu2}.

 \section{ Proof of Theorem 1.1}

     We need the following simple lemma which is a vector-valued function version of Lemma 8.5.5 in \cite {Lu1}.
  \begin{lem} \label{lem 4.1} \ \
     Let $0< \alpha < 1$ and  $\Omega$ be an open set in $\mathbb{R}^n$ with uniformly $C^{r+\alpha}$ boundary, and let $t_0 < t_1$.
     Let $u_i \in C^{r+\alpha, (r+\alpha)/r}(\overline{\Omega} \times [t_0, t_1], \mathbb{R}^l)$ ($i=1,2,\cdot\cdot\cdot$) be a sequence of  $\mathbb{R}^l$-valued functions with
     \begin{equation}
      ||u_i||_{ C^{r+\alpha, (r+\alpha)/r}(\overline{\Omega} \times [t_0, t_1], \mathbb{R}^l) } \leq C,
      \end{equation}
      where the constant $C$ is independent of $i$.   Assume that the sequence $u_i$ converges to $u$ in $C^0( \overline{\Omega} \times [t_0, t_1], \mathbb{R}^l) $.
      Then   $u \in C^{r+\alpha, (r+\alpha)/r}(\overline{\Omega} \times [t_0, t_1], \mathbb{R}^l)$  with
       $||u||_{ C^{r+\alpha, (r+\alpha)/r}(\overline{\Omega} \times [t_0, t_1], \mathbb{R}^l)} \leq C$.
               \end{lem}
      \noindent {\bf Proof}.  $\forall$ $ \Omega' \subset \subset \Omega $ with $C^{r+ \alpha}$ boundary, and $\forall$ $\beta$ with $ 0< \beta < \alpha$,
       by the assumption (4.1) and a consequence of Arzela-Ascoli theorem (compare for example pp. 338-339 in \cite{Al}, Lemma 6.36 in \cite{GT}, and Theorem 2.7 in \cite{Sz} for the elliptic case; the parabolic case is similar), there is a subsequence  of $\{u_i\}$  converges to  $\tilde{u}$ in
     $C^{r+\beta, (r+\beta)/r}(\overline{\Omega'} \times [t_0, t_1], \mathbb{R}^l)$, and
     \begin{equation*}
     ||\tilde{u}||_{C^{r+\beta, (r+\beta)/r}(\overline{\Omega'} \times [t_0, t_1], \mathbb{R}^l)} \leq C.
       \end{equation*}
       By assumption $u_i$ converges to $u$ in $C^0( \overline{\Omega} \times [t_0, t_1], \mathbb{R}^l) $,  then we get that $ \tilde{u}=u$ on  $\overline{\Omega'} \times [t_0, t_1]$ by the uniqueness of  the limit,
    so  $u \in C^{r+\beta, (r+\beta)/r}(\overline{\Omega'} \times [t_0, t_1], \mathbb{R}^l)$ and
      \begin{equation*}
    ||u||_{C^{r+\beta, (r+\beta)/r}(\overline{\Omega'} \times [t_0, t_1], \mathbb{R}^l)} \leq C.
       \end{equation*}
     By the arbitrariness  of $\Omega'$ and $\beta$ the desired result follows.     \hfill{$\Box$}

\vspace*{0.4cm}

Now we prove Theorem 1.1.   Let $M$ be a closed (Riemannian) manifold, $E$ be a (Euclidean) vector bundle over $M$, and $P_t: \text{Dom}(P_t)  \cap  C^\infty(E) \rightarrow C^\infty(E)$, $t\in [0, T]$,
be a  smooth family of smooth partial differential operators of even order $r$, where $\text{Dom}(P_t)$ is an open subset of $C^r(E)$ on which $P_t$ can act. Suppose $u_0\in \text{Dom}(P_0) \cap  C^\infty(E)$ and $P_0$ is strongly  elliptic at $u_0$.  By the compactness of $M$ and smoothness of $P_t$, there are positive constants $\tilde{\delta} (\leq T)$, $\tilde{R}$ and $\tilde{\lambda}$ such that for any $\tilde{u} \in C^r(E)$ with
\begin{equation*}
||\tilde{u}-u_0||_{C^0}+||\nabla \tilde{u}-\nabla u_0||_{C^0}+\cdot\cdot\cdot+||\nabla^r \tilde{u}-\nabla^r  u_0||_{C^0} \leq \tilde{R},
\end{equation*}
we have $\tilde{u} \in \text{Dom}(P_t)$, and
\begin{equation*}
(-1)^{\frac{r}{2}-1} h(\sigma(P_{t*|{\tilde{u}}})(x,\xi)v_x, v_x) \geq \tilde{\lambda} |\xi|^{r}|v_x|^2
\end{equation*}
for any $t\in[0,\tilde{\delta}]$, $(x, \xi)\in T^*M$ and $v_x \in E_x$.
 Given $0< \alpha < 1$, by Theorem 1.2, there exists $\delta >0$, such that the Cauchy problem (1.3) with initial value $u_0$ has  a  unique solution
 $u \in C^{r+\alpha,(r+\alpha)/r}(E_\delta)$.  We may assume that $\delta \leq \tilde{\delta}$. We can assume also that
 \begin{equation*}
||u(\cdot,t)-u_0(\cdot)||_{C^0}+||\nabla u(\cdot,t)-\nabla u_0(\cdot)||_{C^0}+\cdot\cdot\cdot+||\nabla^r u(\cdot,t)-\nabla^r  u_0(\cdot)||_{C^0} \leq \frac{\tilde{R}}{2}
\end{equation*}
 for any $t\in [0,\delta]$ by choosing $\delta$ smaller if necessary.  Then we want to show that in fact $u \in C^\infty(E_\delta)$.
 We need the following local result.
  \begin{prop} \label{prop 4.2} \ \   Suppose that $\Omega$ is a  bounded domain (say an open ball) in $\mathbb{R}^n$ with smooth boundary.  Let $r$ be a positive even number, $k\geq 1$ be an integer, $0< \alpha < 1$ and  $\delta>0$. Assume that $u_0 \in C^{r+k+\alpha} (\overline{\Omega}, \mathbb{R}^l)$,
 and the range of $(u_0, \nabla u_0, \cdot\cdot\cdot, \nabla^r u_0)$  is contained in an open ball
$B(\bar{z}, \frac{\hat{R}}{2})$ in $\mathbb{R}^l \times \mathbb{R}^{nl} \times \cdot\cdot\cdot \times \mathbb{R}^{n^rl}$ for some $\bar{z} \in \mathbb{R}^l \times \mathbb{R}^{nl} \times \cdot\cdot\cdot \times \mathbb{R}^{n^rl}$ and $\hat{R} >0$. Write $Z=B(\bar{z}, \hat{R})$. Furthermore assume that $F$ is a $C^{k+2}$-map from  $\overline{\Omega} \times  [0, \delta] \times Z$ to $\mathbb{R}^l$,  and that  there is a constant $\lambda >0$ such that
    \begin{equation}
  (-1)^{\frac{r}{2}-1}  \frac{\partial F^a}{\partial q^b_{i_1 \cdot\cdot\cdot i_r}}(x,t, z)\xi_{i_1}\cdot\cdot\cdot \xi_{i_r}v^av^b \geq \lambda |\xi|^r|v|^2,
       \end{equation}
  where $(x, t, z) \in  \overline{\Omega} \times  [0, \delta] \times Z$,  $q^b_{i_1 \cdot\cdot\cdot i_r}$ are the coordinates of the $\mathbb{R}^{n^rl}$-component of $z$,  $\xi=(\xi_1,\cdot\cdot\cdot, \xi_n) \in \mathbb{R}^n$, and $v=(v^1,\cdot\cdot\cdot,v^l) \in \mathbb{R}^l$. Let $u \in C^{r+\alpha, (r+\alpha)/r}(\overline{\Omega}\times [0,\delta], \mathbb{R}^l)$ be a solution to the Cauchy problem
  \begin{equation}
\frac{\partial u}{\partial t}=F(x,t, u, \nabla u, \cdot\cdot\cdot, \nabla^ru) \hspace*{6mm} \text{on}  \hspace*{2mm}   \overline{\Omega} \times [0,\delta],  \hspace*{6mm} u(x,0)=u_0(x)  \hspace*{6mm} \text{on}  \hspace*{2mm} \overline{\Omega}
\end{equation}
  with
  \begin{equation}
||u-u_0||_{C^0}+||\nabla u-\nabla u_0||_{C^0}+\cdot\cdot\cdot+||\nabla^r u-\nabla^r  u_0||_{C^0} \leq \frac{\hat{R}}{2}.
\end{equation}
     Then for any domain $\Omega' \subset \mathbb{R}^n$ with $\overline{\Omega'} \subset \Omega$ and $\partial \Omega'$ smooth,  $\nabla^p u \in C^{r+\alpha, (r+\alpha)/r}(\overline{\Omega'} \times [0,\delta], \mathbb{R}^{n^pl})$ for any $p \leq k$.
  \end{prop}
 \noindent {\bf Proof}.  We  use difference quotients and a  bootstrap procedure; compare the proof of Proposition 8.5.6 in \cite {Lu1}, Lemma 14.11 in \cite {Li}, Theorem 8.12.1 in \cite {K},  Theorem 2.5.9 in \cite{Ge}, and Theorem 6.17 in \cite{GT}. We do induction on $k$.
   First consider the case of $k=1$.  Let the hypothesis of the proposition in the case of $k=1$ holds.  Thus $u_0 \in C^{r+1+\alpha} (\overline{\Omega}, \mathbb{R}^l)$,
   and $u \in C^{r+\alpha, (r+\alpha)/r}(\overline{\Omega}\times [0,\delta], \mathbb{R}^l)$ is a solution to (4.3) such that (4.4) holds. Fix any domain $\Omega' \subset \mathbb{R}^n$ with $\overline{\Omega'} \subset \Omega$ and $\partial \Omega'$ smooth. For any $x_0 \in \overline{\Omega'}$, choose $R>0$ such that $B(x_0, 2R) \subset \Omega$. Fix $i$ ($1\leq i \leq n$),  for  $h \in \mathbb{R} \setminus \{0\}$ (not to be confused with the fiber metric $h$ on the bundle $E$) with $|h|< R$, set
  \begin{equation*}
 u_h(x,t)=\frac{u(x+he_i,t)-u(x,t)}{h},  \hspace*{6mm} (x,t) \in B(x_0, R) \times [0,\delta],
   \end{equation*}
  where $e_i$ is the vector in $\mathbb{R}^n$ with the $i$-th component 1 and all the other components 0.  Then $u_h$ satisfies
 $$\begin{array}{l}
 \frac{\partial u_h^a}{\partial t}=A_b^{aI}\partial_Iu_h^b+F_h^a,   \hspace*{6mm} (x,t) \in B(x_0, R)\times [0,\delta],  \\
 u_h(x,0)=\frac{u_0(x+he_i)-u_0(x)}{h}, \hspace*{6mm} x \in B(x_0, R),
  \end{array}$$
 where
 \begin{equation*}
 A_b^{aI}(x,t;h)=\int_0^1   \frac{\partial F^a}{\partial q_I^b}(x+\sigma h e_i, t, \xi_\sigma(x,t;h))d\sigma,
 \end{equation*}

 \begin{equation*}
 F^a_h(x,t)=\int_0^1 \frac{\partial F^a}{\partial x^i}(x+\sigma h e_i, t, \xi_\sigma(x,t;h))d\sigma,
 \end{equation*}
and
  \begin{equation*}
 \xi_\sigma(x,t;h)=\sigma (u,\nabla u,\cdot\cdot\cdot,\nabla^ru)(x+he_i,t)+(1-\sigma)( u,\nabla u,\cdot\cdot\cdot,\nabla^ru)(x,t).
  \end{equation*}
Note that by our assumptions, with $|h|< R$, $(x+\sigma h e_i, \xi_\sigma(x,t;h)) \in  \overline{\Omega} \times Z$ for $(x,t) \in  B(x_0, R) \times  [0, \delta]$ and $\sigma \in [0,1]$.

  Let  $\theta$ be a smooth function on $\mathbb{R}^n$ with
   \begin{equation*}
    \theta \equiv 1 \hspace*{2mm} \text {in} \hspace*{2mm}  B(x_0, \frac{R}{2}),  \hspace*{2mm} \theta \equiv 0 \hspace*{2mm}  \text{outside}  \hspace*{2mm}  B(x_0,R).
   \end{equation*}
Define a family of $\mathbb{R}^l$-valued functions $\tilde{v}$ on $\overline{\Omega} \times [0,\delta]$ with parameter $h$ by
  \begin{equation*}
  \tilde{v}(x,t;h)=\theta (x)u_h(x,t)  \hspace*{2mm} \text {in} \hspace*{2mm} B(x_0, R)\times [0,\delta], \hspace*{2mm}  \text {and}  \hspace*{2mm}      \tilde{v}(x,t;h)=0 \hspace*{2mm} \text{elsewhere}.
     \end{equation*}
   Then for each $h\neq 0$ with $|h|< R$, $\tilde{v}$ satisfies the system
    $$\begin{array}{l}
  \frac{\partial \tilde{v}}{\partial t}=L_t(\tilde{v})+\tilde{F},   \hspace*{6mm}   (x,t) \in  \overline{\Omega}\times [0,\delta],  \\
  \tilde{v}(x,0;h)=\tilde{v}_0(x;h),  \hspace*{6mm}  x \in \overline{\Omega },    \\
  \tilde{v}(x,t;h)=0,   \hspace*{6mm}  (x,t) \in \partial \Omega \times [0,\delta],

  \end{array}$$
   where
   \begin{equation*}
   (L_t\tilde{v})^a= A_b^{aI}(x,t;h)\partial_I\tilde{v}^b,
  \end{equation*}
  \begin{equation*}
    \tilde{F}=\theta F_h+(\theta L_t(u_h)-L_t(\theta u_h)),
   \end{equation*}
   and
    \begin{equation*}
      \tilde{v}_0(x;h)=\theta(x)\frac{u_0(x+he_i)-u_0(x)}{h} \hspace*{2mm}  \text{in}    \hspace*{2mm} B(x_0,R),    \hspace*{2mm} \text{ and}   \hspace*{2mm}
      \tilde{v}_0(x;h)  \hspace*{2mm} \text{is} \hspace*{2mm} 0 \hspace*{2mm}  \text{elsewhere}.
          \end{equation*}

  Note that due to our assumptions the operators $L_t$ ($t\in [0, \delta]$) are strongly elliptic with the ellipticity constant independent of $t$ and $h$. Observe that any $\nabla^pu_h$ appearing in the  difference
   $\theta L_t(u_h)-L_t(\theta u_h)$ must have $p\leq r-1$.
  By assumption  $u \in C^{r+\alpha, (r+\alpha)/r}(\overline{\Omega} \times [0,\delta], \mathbb{R}^l)$, so
  the $C^{\alpha,\alpha/r}(\overline{\Omega} \times [0,\delta], \mathbb{R}^{n^pl})$-norms of $\nabla^p u_h$  ($p\leq r-1$)   are bounded by a constant independent of $h$. It follows that the $C^{\alpha,\alpha/r}(\overline{\Omega} \times [0,\delta], \mathbb{R}^l)$-norm of $\tilde{F}$  is bounded by a constant independent of $h$. Also note that
  the $C^{r+\alpha, (r+\alpha)/r}(\overline{\Omega})$-norm of
   $\tilde{v}_0(\cdot;h)$  is bounded by a constant independent of $h$. Then it follows from Theorem 4.9   of \cite {So} (see also Theorem 10.1 in Chapter VII of \cite{LSU}) that
    the $C^{r+\alpha, (r+\alpha)/r}(\overline{\Omega} \times [0,\delta], \mathbb{R}^l)$-norm of $\tilde{v}$ is bounded by a constant $C$ independent of $h$.
     (Note that $\tilde{v}_0$ is 0 in a neighborhood of $\partial \Omega $ in $\overline{ \Omega} $,   so the initial-boundary value problem
     for $\tilde{v}$ above satisfies the required complementarity and compatibility conditions.)  Clearly $u_h \rightarrow \frac{\partial u}{\partial x^i} $ in $C^0(B(x_0,R) \times [0, \delta], \mathbb{R}^l)$ as $h\rightarrow 0$.
   Then by Lemma 4.1 we have $\frac{\partial u}{\partial x^i} \in C^{r+\alpha, (r+\alpha)/r}(B(x_0, \frac{R}{2}) \times [0,\delta], \mathbb{R}^l)$, and $||\frac{\partial u}{\partial x^i} ||_{C^{r+\alpha, (r+\alpha)/r}(B(x_0, R/2) \times [0,\delta],\mathbb{R}^l)} \leq C$. So  the proposition is true for the case of $k=1$.

 Suppose the proposition  is true for the case of $k=j$. Now consider the case of $k=j+1$.  Let the hypothesis of the proposition in the case of $k=j+1$ holds.  In particular, $u_0 \in C^{r+j+1+\alpha} (\overline{\Omega}, \mathbb{R}^l)$,
   and $u \in C^{r+\alpha, (r+\alpha)/r}(\overline{\Omega}\times [0,\delta], \mathbb{R}^l)$ is a solution to (4.3)  such that (4.4) holds. Fix  a domain $\Omega' \subset \mathbb{R}^n$ with $\overline{\Omega'} \subset \Omega$.  By induction assumption, for any domain $\Omega_1 \subset \mathbb{R}^n$ with $\overline{\Omega'} \subset  \Omega_1 \subset \overline{\Omega_1} \subset \Omega$  and $\partial \Omega_1$ smooth,  $\nabla^p u \in C^{r+\alpha, (r+\alpha)/r}(\overline{\Omega_1} \times [0,\delta], \mathbb{R}^{n^pl})$ for any $p \leq j$.

Fix  $1 \leq i_1, \cdot\cdot\cdot, i_j \leq n$. Let $v=\frac{\partial ^{j}u}{\partial x^{i_1} \cdot\cdot\cdot \partial x^{i_j}}$.
 We differentiate (4.3) w.r.t. the variables $x^{i_1}, \cdot\cdot\cdot, x^{i_j}$ successively and get
$$\begin{array}{l}
\frac{\partial v^a}{\partial t}(x,t)=\frac{\partial F^a}{\partial q_{I}^b}\partial_{I}v^b(x,t)+\hat{F}^a,
 \hspace*{6mm} (x,t) \in \overline{\Omega_1} \times [0,\delta],  \\
v(x,0)=\frac{\partial ^{j}u_0}{\partial x^{i_1} \cdot\cdot\cdot \partial x^{i_j}}(x), \hspace*{6mm} x \in \overline{\Omega_1}.
 \end{array}$$
Here,   $I$ is any  multi-index with $|I|\leq r$, each coefficient $\frac{\partial F^a}{\partial q_{I}^b}$  is evaluated at the point $(x,t,u(x,t), \nabla u(x,t), \cdot\cdot\cdot, \nabla^ru(x,t))$, and any (spatial) derivative of $u$ contained in $\hat{F}^a$ has order $\leq r+j-1$. (Of course here as usual the repeated indices mean summation.) One can easily check that each  $\frac{\partial F^a}{\partial q_{I}^b}$ and $\hat{F}^a$ are at least of the class  $C^{1+\alpha, (1+\alpha)/r}$. Moreover the initial value $v(\cdot,0)$ is of the class $C^{r+1+\alpha}$. By the same argument via difference quotients as in the case of $k=1$ we get that $||\frac{\partial v}{\partial x^i} ||_{C^{r+\alpha, (r+\alpha)/r}(\overline{\Omega'} \times [0,\delta],\mathbb{R}^l)} \leq C$ for any $1\leq i \leq n$.
Then the desired result in the case of $k=j+1$  follows.
  \hfill{$\Box$}

\vspace*{0.4cm}

Note that under the hypothesis of Proposition 4.2, we can also obtain higher regularity of $u$ w.r.t. the time variable $t$ by using the difference quotient trick for $t$, Proposition 4.2, the system (4.3) and the systems derived by differentiating  (4.3) w.r.t. $t$ (many times) and Theorem 4.9   of \cite {So}; similarly we can treat the space-time mixed derivatives. (Consult the proof of Theorem  2.5.9 in \cite{Ge}.) Then the proof of Theorem 1.1 is finished. \hfill{$\Box$}


\hspace *{0.4cm}

\vspace*{0.4cm}

Laboratory of Mathematics and Complex Systems (Ministry of Education),

School of Mathematical Sciences, Beijing Normal University,

Beijing 100875,  People's Republic of China

 E-mail address: hhuang@bnu.edu.cn

\end{document}